\documentclass{article}
\usepackage{amsmath,amsfonts,amssymb,epsf,chicago,pb-diagram,lamsarrow,pb-lams,theorem}

\newcommand{\rf}[1]{(\expandafter\ref{#1})}
\newcommand{\ct}[1]{\citeANP{#1}~[\citeyearNP{#1}]}

\newcommand{\lb}[1]{\expandafter\label{#1}}

\newcommand\mfk\mathfrak
\newcommand\mcl\mathcal
\newcommand\mbb\mathbb
\newcommand\onm\operatorname
\newcommand{\Id}{{\mathchoice{\mbox{Id}}{\mbox{Id}}{\mbox{\scriptsize Id}}{\mbox{\scriptsize Id}}}}
\newcommand{\set}[2]{\bigl\{\,#1:#2\,\bigr\}}
\newcommand{\bset}[2]{\bigl\{\,#1:#2\,\bigr\}}

\newtheorem{theorem}{Theorem}
\newtheorem{corollary}{Corollary}
\newtheorem{definition}{Definition}
\newtheorem{proposition}{Proposition}
\newtheorem{lemma}{Lemma}
\theorembodyfont{\rmfamily}
\newtheorem{proof}{Proof}
\newcommand{\qed}{\\\mbox{}\hfill$\square$}

\begin{document}\allowdisplaybreaks
\title{The Transversal Relative Equilibria of a Hamiltonian System 
with Symmetry}
\author{G. W. Patrick\\
Department of Mathematics and Statistics\\
University of Saskatchewan\\Saskatoon, Saskatchewan, 
S7N~5E6\\Canada
\and
R. M. Roberts\\
Mathematics Institute\\
University of Warwick\\
Coventry, CV4 7AL\\
United Kingdom}
\date{2 February 1999}
\maketitle

\begin{abstract}
We show that, given a certain transversality condition, the set of
relative equilibria $\mcl E$ near $p_e\in\mcl E$ of a Hamiltonian
system with symmetry is locally Whitney-stratified by the conjugacy
classes of the isotropy subgroups (under the product of the coadjoint
and adjoint actions) of the momentum-generator pairs $(\mu,\xi)$ of
the relative equilibria. The dimension of the stratum of the conjugacy
class $(K)$ is $\dim G+2\dim Z(K)-\dim K$, where $Z(K)$ is the
center of $K$, and transverse to this stratum $\mcl E$ is locally
diffeomorphic to the commuting pairs of the Lie algebra of $K/Z(K)$.
The stratum $\mcl E_{(K)}$ is a symplectic submanifold of $P$ near
$p_e\in\mcl E$ if and only if $p_e$ is nondegenerate and $K$ is a
maximal torus of $G$. We also show that there is a dense subset
of $G$-invariant Hamiltonians on $P$ for which all the relative
equilibria are transversal. Thus, generically, the types of
singularities that can be found in the set of relative equilibria of a
Hamiltonian system with symmetry are those types found amongst the
singularities at zero of the sets of commuting pairs of certain Lie
subalgebras of the symmetry group.
\end{abstract}

\section*{Introduction} 

Let $G$ be a connected, compact Lie group, with Lie algebra $\mfk g$.
 The group $G$ acts via its adjoint representation on $\mfk g$ and by
 the dual coadjoint representation on the dual $\mfk g^*$ of $\mfk
 g$. Let $G$ act freely and symplectically on a symplectic phase space
 $(P,\omega)$ endowed with an equivariant momentum mapping
 $J:P\rightarrow\mfk g^*$. Since $G$ is acting freely, every point in
 $P$ is a regular point of $J$. Fix a $G$-invariant Hamiltonian $H$ on
 $P$. A relative equilibrium of $H$ with generator $\xi_e\in\mfk g$ is
 a point $p_e\in P$ for which the evolution under the flow generated
 by $H$ is $t\mapsto \onm{exp}(t\xi_e).p_e$. Let $G_{\mu_e}$ and
 $G_{\xi_e}$ respectively denote the isotropy subgroups of
 $\mu_e=J(p_e)\in\mfk g^*$ and $\xi_e\in\mfk g$ under the coadjoint
 and adjoint actions of $G$. \ct{PatrickGW-1995.1} has proved that if
 $G_{\mu_e}\cap G_{\xi_e}$ is a maximal torus of $G$ and a
 nondegeneracy condition holds then the set of relative equilibria
 near $p_e$ is a manifold of dimension $\onm{dim}G+\onm{rank}G$. This
 manifold is symplectic under an
 additional nonresonance condition and the condition that $\mu_e$ is
 generic.

Recently \ct{OrtegaJPRatiuTS-1997.2} have applied these results,
combined with singular reduction techniques, to obtain
manifolds of relative equilibria even when $p_e$ is not regular. In
an interesting article with a number of new ideas,
\ct{MontaldiJA-1997.1} proves that, regardless of whether or not $G$
acts freely, certain extremal relative equilibria must persist to 
relative equilibria on symplectic reduced
spaces near to the one through $p_e$, and exactly counts and classifies
the relative equilibria near $p_e$ when $G$ does act freely and the
generator $\xi_e$ is regular. Other results on the existence of manifolds
of relative equilibria through a given non-regular relative equilibrium
have been obtained by \ct{RobertsRMSousaDiasMER-1997.1} and
\ct{LermanESingerSF-1998.1}.

We define the {\it type} of a pair $(\xi,\mu)\in\mfk g^*\oplus\mfk g$
to be the conjugacy class in $G$ of $G_{\mu}\cap G_{\xi}$ and the
type of a relative equilibrium $p_e$ to be the type of its
momentum-generator pair $(\mu_e,\xi_e)$. Generally we denote the
conjugacy class of a subgroup $K\subset G$ by $(K)$, and so we say that
the type of $p_e$ is $(K_e)$, where $K_e=G_{\mu_e}\cap G_{\xi_e}$. The
main result of this paper is a description of the complete set of
relative equilibria near $p_e$ in the relaxed situation where the type
of $p_e$ is not necessarily a maximal torus of $G$, but where $p_e$ is
still regular. One of our opinions, indicated first
by~\ct{PatrickGW-1995.1}, and born out by our work, is that the
momentum-generator pair is a more important object than either the
momentum or generator taken singly.

Consider the case of a single `generic' rigid body with distinct
principle moments of inertia $I_1$, $I_2$ and $I_3$. This is an
$SO(3)$-invariant system. The relative equilibria consist of steady
rotations about the three principle axes of inertia. After the $SO(3)$
symmetry is removed by reduction to the Euler equations, the set of
relative equilibria is the union of the three coordinate axes of $\mbb
R^3$. The origin has type $SO(3)$ and is the sole bifurcation point of
the set of relative equilibria. The remaining equilibria have
type~$\bigl(SO(2)\bigr)$, the `smallest' possible. On the evidence of
this example, one might expect general $SO(3)$-invariant systems to
have a zero dimensional set of relative equilibria of type $SO(3)$
from which bifurcate three distinct branches of relative equilibria of
type $\bigl(SO(2)\bigr)$. {\em However, the results we obtain imply
this example is a poor guide to the general theory because it is
nongeneric. Generically, in the class of Hamiltonian systems with
symmetry group $SO(3)$, there are no regular relative equilibria of
type~$SO(3)$. In particular, this implies that there are no regular
equilibria.}

More generally let $\mcl E\subset P$ denote the set of relative
equilibria of a generic $G$-invariant Hamiltonian system on $P$ and
let $\mcl E/G\subset P/G$ denote its quotient by the action of
$G$. Then we show that both $\mcl E$ and $\mcl E/G$ are stratified by
the types of the relative equilibria and that the quotient of a
non-empty stratum consisting of relative equilibria of type $(K)$ has
dimension equal to $2\dim Z(K)-\dim K$, where $Z(K)$ is the center of
$K$. Moreover transverse to this
stratum the set $\mcl E/G$ is locally isomorphic to the set of
commuting pairs of the Lie algebra of $K/Z(K)$. See Theorems~\ref{mainth1} and~\ref{mainth2} for precise statements.

The types of transversal relative equilibria are constrained by
the condition $2\dim Z(K)\ge\dim K$, a condition which fails for 
$K=SO(3)$ in the rigid body example.
For another example consider a generic rigid body with
three rotors aligned with its (three orthogonal) principle axes. This
has symmetry group $SO(3)\times SO(2)^3$. There is, up to symmetry,
exactly one relative equilibrium of type $SO(3)\times SO(2)^3$, namely
the equilibrium where neither the body nor the rotors move. At this
equilibrium the set of relative equilibria has a singularity,
bifurcating from which are very many relative equilibria of (generic)
type $\bigl(SO(2)^4\bigr)$. {\em Theorems~\ref{mainth1}, \ref{mainth2}
and~\ref{mainth3} imply that, for generic Hamiltonian systems with
symmetry group $SO(3)\times SO(2)^3$, the relative equilibria of type
$SO(3)\times SO(2)^3$ form isolated group orbits which are
singularities of $\mcl E/G$ and that near each of these singularities
$\mcl E/G$ is diffeomorphic to the set of commuting pairs of the Lie
algebra $so(3)$.}

Here is a summary of this work: Let $(\mfk g^*\oplus\mfk
g)^c\subseteq\mfk g^*\oplus\mfk g$ be the set of pairs $(\mu,\xi)$
such that $\onm{coad}_\xi\mu = 0$, where `$\onm{coad}$' denotes the
coadjoint representation of $\mfk g$ on $\mfk g^*$. Any $G$-invariant
inner product on $\mfk g$ defines an isomorphism between $\mfk g$ and
$\mfk g^*$ and between the set of commuting pairs of $\mfk g$ and the
subset $(\mfk g^*\oplus\mfk g)^c$. In Section~\ref{g*g} we describe a
stratification of $(\mfk g^*\oplus\mfk g)^c$ and in
Section~\ref{transRE} we use this to construct stratifications of
certain subsets $\mcl T^c$ of the tangent bundle and $\mcl K^{oc}$ of
the cotangent bundle of $P$. A relative equilibrium is defined to be
{\em transversal} if at that point the vector field is transversal to
the stratification of $\mcl T^c$ or, equivalently, the derivative of
the Hamiltonian is transverse to the stratification of $\mcl K^{oc}$.
Some more or less standard results from transversality and
stratification theory are then used to deduce Theorems~\ref{mainth1},
\ref{mainth2} and~\ref{mainth3}. In Section~\ref{linear} we describe
a normal form for the linearization of a vector field at a transversal
relative equilibrium and use this to show that near a transversal
relative equilibrium $p_e$ the submanifold of $P$ consisting of
relative equilibria of type $(K_e)$ is symplectic if and only if $p_e$
is nondegenerate and $K_e$ is a maximal torus in $G$
(Theorem~\ref{mainth4}). This is a generalization and partial converse
of a result of \ct{PatrickGW-1995.1}.

\section{The Stratified Structure of $(\mfk g^*\oplus\mfk g)^c$}\lb{g*g}

The momentum-generator pair $(\mu_e,\xi_e)\in\mfk g^*\oplus\mfk g$ of a
relative equilibrium of a $G$-invariant Hamiltonian system satisfies the
relationship $\onm{coad}_{\xi_e}\mu_e=0$. In this section we describe the
structure of the set of all pairs satisfying this relationship:
\begin{equation*}
(\mfk g^*\oplus\mfk g)^c=\set{(\mu,\xi)\in\mfk g^*\oplus\mfk g}
{\onm{coad}_\xi\mu = 0}.
\end{equation*}

Note that $(\mfk g^*\oplus
\mfk g)^c$ is a $G$-invariant subvariety of $\mfk g^*\oplus\mfk g$ under the
product of the coadjoint and adjoint actions.
For any $(\mu,\xi)\in\mfk g^*\oplus\mfk g$ let $G_{\mu,\xi}$ denote the
isotropy subgroup of $(\mu,\xi)$ for the product of the coadjoint and
adjoint actions. Clearly $G_{\mu,\xi} = G_\mu\cap G_\xi$ where
$G_\mu$ is the isotropy subgroup of $\mu$ for the coadjoint action of
$G$ and $G_\xi$ is the isotropy subgroup of $\xi$ for the adjoint
action. Let $\mfk g_\mu$ and $\mfk g_\xi$ denote the Lie algebras of $G_\mu$
and $G_\xi$, respectively, and $\mfk g_{\mu}^{*}$ and $\mfk g_{\xi}^{*}$
their dual spaces.

\begin{lemma}\lb{g*gisgs}\mbox{}
\begin{enumerate}
\item 
$(\mu,\xi)\in(\mfk g^*\oplus\mfk g)^c$ if and only if
$\xi\in\mfk g_\mu\subset\mfk g$, and $G_{\mu,\xi}=(G_\mu)_\xi$, the isotropy
subgroup of $\xi$ for the adjoint action of $G_\mu$ on $\mfk g_\mu$.
\item 
A point $(\mu,\xi)\in\mfk g^*\oplus\mfk g$ lies in $(\mfk
g^*\oplus\mfk g)^c$ if and only if $G_{\mu,\xi}$ contains a maximal
torus of $G$.
\item 
If $(\mu,\xi)\in\mfk (g^*\oplus\mfk g)^c$ then $\mfk g_{\mu,\xi}=\mfk
g_\mu\cap\mfk g_\xi$ is abelian if and only if $G_{\mu,\xi}$ is a
maximal torus of $G$.
\end{enumerate}
\end{lemma}

\begin{proof}
We have
\begin{equation*}
(\mu,\xi)\in(\mfk g^*\oplus\mfk g)^c\ \Leftrightarrow\ 
\onm{coad}_\xi\mu= 0\ \Leftrightarrow\ \xi\in\mfk g_\mu.
\end{equation*}
Moreover
\begin{equation*}
G_{\mu,\xi}=G_\mu\cap G_\xi=
\set{g\in G_\mu}{\onm{Ad}_g\xi=\xi}
\end{equation*}
and part~1 follows from the fact that the adjoint action of $G$ on
$\mfk g$ restricts to the adjoint action of $G_\mu$ on $\mfk g_\mu$.

For part~2, suppose first that $(\mu,\xi)\in(\mfk g^*\oplus\mfk g)^c$.
Then, by part 1, $G_{\mu,\xi}$ is an isotropy subgroup of the adjoint
action of $G_\mu$ and so contains a maximal torus of $G_\mu$. However
$G_\mu$ is an isotropy subgroup of the coadjoint action of $G$ and so
maximal tori of $G_\mu$ are also maximal tori of $G$.

Conversely, suppose $G_{\mu,\xi}$ contains a maximal torus $T$ of $G$.
Then $(\mu,\xi)\in\onm{fix}(T;\mfk g^*\oplus\mfk g)=\mfk t^*\oplus\mfk
t$ where $\mfk t$ is the Lie algebra of $T$. Since $T$ is abelian it
follows that $\onm{coad}_\xi\mu= 0$ and so $(\mu,\xi)\in(\mfk
g^*\oplus\mfk g)^c$.

For part~3, if $G_{\mu,\xi}$ is a maximal torus then its Lie algebra
$\mfk g_{\mu,\xi}$ is abelian. Conversely, suppose $\mfk g_{\mu,\xi}$
is abelian. Then by Proposition 4.25 of~\ct{AdamsJF-1969.1} and by
part~1, $G_{\mu,\xi}$ is connected and so is also abelian. Thus
$G_{\mu,\xi}$ is a torus that contains a maximal torus, and so is a
maximal torus.\qed
\end{proof}

For any set $X$ with an action of $G$ we denote the subset consisting
of points with isotropy subgroup conjugate to $K\subset G$ by
$X_{(K)}$. We will denote the set of conjugacy classes of isotropy subgroups of
the action of $G$ on $(\mfk g^*\oplus\mfk g)^c$ by $\mcl I^c$.
By~Lemma~\ref{g*gisgs}, $\mcl I^c$ is the set of conjugacy classes of
isotropy subgroups of the action of $G$ on $\mfk g^*\oplus\mfk g$
which contain a maximal torus of $G$, and
\begin{equation*}
(\mfk g^*\oplus\mfk g)^c =\bigsqcup_{(K)\in\mcl I^c}
 (\mfk g^*\oplus\mfk g)_{(K)}.
\end{equation*}
By the general theory of actions of a compact group, each
set $(\mfk g^*\oplus\mfk g)_{(K)}$ is a manifold and the collection of
these manifolds is a Whitney regular stratification of $\mfk
g^*\oplus\mfk g$. Since $(\mfk g^*\oplus\mfk g)^c$ is the union of
all the strata of $\mfk g^*\oplus\mfk g$ which are smaller than the
stratum of the maximal tori, it follows that $\bigl\{(\mfk
g^*\oplus\mfk g)_{(K)}:(K)\in\mcl I^c\bigr\}$ is a Whitney regular
stratification of $(\mfk g^*\oplus\mfk g)^c$. The general theory also
tells us that the quotients of the strata by the action of $G$ are
smooth manifolds and that the set of these strata is a Whitney regular
stratification of the orbit space $(\mfk g^*\oplus\mfk g)/G$. It
follows that the set of quotients of the strata $(\mfk g^*\oplus\mfk
g)_{(K)}$ with $(K)\in\mcl I^c$ is a Whitney regular stratification of
$(\mfk g^*\oplus\mfk g)^c/G$.

The next proposition describes the local structure of these
stratifications. We fix a $G$-invariant inner product on $\mfk g$,
which determines a $G$-equivariant isomorphism between $\mfk g$ and
$\mfk g^*$ and hence also a $G$-invariant inner product on $\mfk g^*$.
The $G$-invariant inner product also defines an isomorphism between
$(\mfk g^*\oplus\mfk g)^c$ and the set of commuting pairs of $\mfk g$:
\begin{equation*}
(\mfk g\oplus\mfk g)^c=\set{(\eta,\xi)\in\mfk g\oplus\mfk g}
{\onm{ad}_\xi\eta = [\xi,\eta] = 0}.
\end{equation*}
If $\mfk k$ is a subspace of $\mfk g$ the inner product on $\mfk g^*$ is
used to identify $\mfk k^*$ with the subspace
$\onm{ann}(\mfk k)^\perp$ of $\mfk g^*$. Here $\onm{ann}(\mfk k)$ is the
annihilator of $\mfk k$ in $\mfk g^*$ and `$\perp$' denotes the orthogonal
complement with respect to the $G$-invariant inner product on $\mfk g^*$.

Suppose $(\mu,\xi)\in (\mfk g^*\oplus\mfk g)^c$
and let $K = G_{\mu}\cap G_{\xi}$. Denote the Lie algebra of $K$ by
$\mfk k$. Let $\chi : U\rightarrow G$ be a section of the natural
projection $G\rightarrow G/K$, defined on an open neighborhood
$U$ of $K$ in $G/K$ and such that $\chi(K)$ is the identity in $G$.
Give $\mfk g^*\oplus\mfk g$ and $(\mfk g^*\oplus\mfk g)^c$ the
stratifications obtained by taking the
connected components of the orbit type strata for the actions of $G$.
Similarly give $\mfk k^*\oplus\mfk k$ and $(\mfk k^*\oplus\mfk k)^c$
the corresponding stratifications obtained
from the actions of $K$. Extend these stratifications to
$U\times (\mfk k^*\oplus\mfk k)$ and
$U\times (\mfk k^*\oplus\mfk k)^c$ by taking the product of each stratum
with $U$. We will always consider any open subset of a stratified set to
be automatically endowed with the stratification obtained by taking the
intersections of the strata with the open subset.

\begin{proposition}\lb{prop1}
Under the above assumptions, there exists a
$G$-invariant open neighborhood $W$ of $(\mu,\xi)$ in
$\mfk g^*\oplus\mfk g$ and a $K$-invariant open neighborhood $V$ of the
origin in $\mfk k^*\oplus\mfk k$ such that the map
\begin{align*}
&\Sigma:U\times V\rightarrow\mfk g^*\oplus\mfk g\\
&\bigl(u,(\nu,\eta)\bigr)\mapsto\chi(u).\bigl((\mu,\xi)+(\nu,\eta)\bigr)
\end{align*}
is an embedding of $U\times V$ into $W$ which restricts to an isomorphism
of smoothly stratified spaces between
$U\times\bigl(V\cap(\mfk k^*\oplus\mfk k)^c\bigr)$ and
$W\cap(\mfk g^*\oplus\mfk g)^c$.
\end{proposition}

\begin{proof}
By the theory of compact group actions there is a $K$-invariant
neighborhood $W$ of $(\mu,\xi)$ in $\mfk g^*\oplus\mfk g$ such
that
\begin{equation*}
S_{\mu,\xi}=W\cap\bigl((\mu,\xi)+(\mfk g.(\mu,\xi))^\perp\bigr)
\end{equation*}
is a $K$-invariant slice to the action of $G$ on $\mfk g^*\oplus\mfk
g$ at $(\mu,\xi)$. A straightforward calculation shows that $\mfk
k^*\oplus\mfk k$ is contained in $(\mfk g.\bigl(\mu,\xi)\bigr)^\perp$. Let
$V$ be the intersection of $\mfk k^*\oplus\mfk k$ with
$S_{\mu,\xi}-(\mu,\xi)$. Then it follows from the slice
theorem that $\Sigma$ is an embedding of $U\times V$ into $W$.

It remains to be proved that 
\begin{equation*}
S_{\mu,\xi}\cap(\mfk g^*\oplus\mfk g)^c=(\mu,\xi)+V
\cap(\mfk k^*\oplus\mfk k)^c.
\end{equation*} 
The righthand side of this equation is clearly contained in the
lefthand side, so it is sufficient to prove that the lefthand side is
contained in the righthand side. Note that for any $(\nu,\eta)\in(\mfk
g^*\oplus\mfk g)$ we always have $\nu\in\mfk g_\nu^*$ and $\eta\in\mfk
g_\eta$. If $(\nu,\eta)\in (\mfk g^*\oplus\mfk g)^c$ then we also have
$\nu\in\mfk g_\eta^*$ and $\eta\in\mfk g_\nu$, and hence $\nu\in\mfk
g_\eta^*\cap\mfk g_\nu^*$ and $\eta\in\mfk g_\eta\cap\mfk g_\nu$. If
$(\nu,\eta)$ lies in the slice $S_{\mu,\xi}$ then the isotropy
subgroup $G_{\nu,\eta}$ must be contained in $K$, and so $\mfk
g_\eta\cap\mfk g_\nu\subset\mfk k$. It follows that if $(\nu,\eta)\in
S_{\mu,\xi}\cap(\mfk g^*\oplus\mfk g)^c$ then $(\nu,\eta)$ must
lie in $\mfk k^*\oplus\mfk k$, and hence in $(\mfk k^*\oplus\mfk
k)^c$, as required.\qed
\end{proof}

This proposition shows that the local structure of the stratification
of $(\mfk g^*\oplus\mfk g)^c$ near a point with isotropy subgroup $K$
can be reduced to that of the set $(\mfk k^*\oplus\mfk k)^c$. For any
$(K)\in\mcl I^c$ let $\mcl I^c(K)$ denote the subset of $\mcl I^c$
consisting of orbit types $(K')$ for which the closure of the stratum
$(\mfk g^*\oplus\mfk g)_{(K')}$ contains $(\mfk g^*\oplus\mfk
g)_{(K)}$. Then the proof of the Proposition implies that $\mcl
I^c(K)$ is the set of conjugacy classes in $G$ of the isotropy
subgroups of the action of $K$ on $(\mfk k^*\oplus\mfk k)^c$.

The next result, which is valid for any connected compact Lie group
$K$, gives a further reduction. Let $Z(K)$ denote the center of $K$.
Denote the quotient $K/Z(K)$ by $L$, its Lie algebra by $\mfk l$ and the
Lie algebra of $Z(K)$ by $\mfk z$. Identify $\mfk z^*$ with the
subspace $\onm{ann}(\mfk z)^\perp$ of $\mfk k^*$.

\begin{proposition}\lb{prop2}
A $K$-invariant inner product on $\mfk k$ determines a
$K$-equivariant linear isomorphism
\begin{equation*}
\sigma:\mfk k^*\oplus\mfk k\cong (\mfk l^*\oplus\mfk l)
\oplus(\mfk z^*\oplus\mfk z)
\end{equation*}
where the action of $K$ on $\mfk k^*\oplus\mfk k$ is the product of
coadjoint and adjoint actions, the action on $\mfk l^*\oplus\mfk l$
factors through the product of the coadjoint and adjoint actions of
$L$ and the action on $\mfk z^*\oplus\mfk z$ is trivial. Moreover
$\sigma$ maps $\onm{fix}(K;\mfk k^*\oplus\mfk k)$ isomorphically to
$\mfk z^*\oplus\mfk z$ and $(\mfk k^*\oplus\mfk k)^c$ to
$(\mfk l^*\oplus\mfk l)^c\oplus(\mfk z^*\oplus\mfk z)$. 
\end{proposition}

\begin{proof}
Since $K$ acts by the adjoint action on $\mfk k$, $\mfk
z=\onm{fix}(K;\mfk k)$. The natural projection from $\mfk k$ to $\mfk
k/\mfk z\cong\mfk l$ is equivariant with respect to the natural
actions of $K$ on $\mfk k$ and $\mfk l$ and so induces an isomorphism
of representations between $\onm{fix}(K;\mfk k)^\perp$ and $\mfk
l$. Hence $\mfk k$ is isomorphic as a representation of $K$ to $\mfk
l\oplus\mfk z$. The invariant inner product translates this into an
isomorphism $\mfk k^*\cong\mfk l^*\oplus\mfk z^*$ and putting the two
together gives $\sigma$. Since the adjoint action of $\mfk z$ on $\mfk
k$ is trivial $\sigma$ maps $(\mfk k^*\oplus\mfk k)^c$ to $(\mfk
l^*\oplus\mfk l)^c\oplus(\mfk z^*\oplus\mfk z)$. \qed\end{proof}

\noindent Thus the orbit type stratification of $(\mfk k^*\oplus\mfk
k)^c$ is isomorphic to the stratification of $(\mfk l^*\oplus\mfk
l)^c\oplus(\mfk z^*\oplus\mfk z)$ obtained by taking the products of
the orbit type strata of $(\mfk l^*\oplus\mfk l)^c$ with $(\mfk
z^*\oplus\mfk z)$.

\begin{corollary}
If $(K)\in\mcl I^c$ then
the dimension of $(\mfk g^*\oplus\mfk g)^c_{(K)}$ is equal to
$\dim G+2\dim Z(K)-\dim K$. 
\end{corollary}

\begin{proof}
By Proposition~\ref{prop1}
\begin{equation*}
\onm{dim}(\mfk g^*\oplus\mfk g)_{(K)}^c=
\onm{dim}U+\onm{dim}(\mfk k^*\oplus\mfk k)_{(K)}^c.
\end{equation*}
The dimension of $U$ is $\onm{dim}G-\onm{dim}K$ while by
Proposition~\ref{prop2} the dimension of
\begin{equation*}
(\mfk k^*\oplus\mfk k)_{(K)}^c=\onm{fix}(K;\mfk k^*\oplus\mfk k)=
\mfk z^*\oplus\mfk z
\end{equation*}
is $2\onm{dim}\mfk z$. This gives the result.
\qed\end{proof}

We note that these structure results also follow from the
results  of~\ct{ArmsJMMarsdenJEMoncriefV-1981.1}
and~\ct{SjamaarRLermanE-1991.1} on the structure of level sets of
momentum maps. Indeed, the product of the coadjoint and adjoint
actions of $G$ on $\mfk g^*\oplus\mfk g$ is symplectic with respect to
the natural `cotangent bundle' symplectic structure, the map
$(\nu,\eta)\mapsto\onm{coad}_\eta\nu$ is an equivariant momentum map for
this action, and $(\mfk g^*\oplus\mfk g)^c$ is its zero level set.
Moreover the symplectic normal space to the group orbit through any
point $(\mu,\xi)\in(\mfk g^*\oplus\mfk g)^c$ with isotropy
subgroup $K$ can be identified with $\mfk k^*\oplus\mfk k$. This
approach also shows that the strata $(\mfk g^*\oplus\mfk g)_{(K)}^c$
are symplectic submanifolds of $\mfk g^*\oplus\mfk g$.

\section{Transversal Relative Equilibria}\lb{transRE}

In this section we first define what we mean by a transversal relative
equilibrium (Definition~\ref{defn}) and then give the main results on
the structure of the space of relative equilibria near a transversal
relative equilibrium (Theorems~\ref{mainth1} and~\ref{mainth2}) and on
the genericity of the set of Hamiltonians for which all the relative
equilibria are transversal (Theorem~\ref{mainth3}).

Define the map $\widetilde J:P\times\mfk g\rightarrow\mfk
g^*\oplus\mfk g$ by $(p,\xi)\mapsto\bigl(J(p),\xi\bigr)$. Since $G$ is
acting freely the map $J$ is a $G$-equivariant submersion and hence so
also is $\widetilde{J}$. Define
\begin{equation*}
(P\times\mfk g)^c=\widetilde J^{-1}(\mfk g^*\oplus\mfk g)^c=
\set{(p,\xi)\in P\times\mfk g}{\onm{coad}_\xi J(p)=0}.
\end{equation*}
If $p_e$ is a relative equilibrium of a $G$-invariant Hamiltonian system
on $P$ with generator $\xi_e$ then $(p_e,\xi_e)\in(P\times\mfk g)^c$. The fact
that $\widetilde{J}$ is a submersion implies that the Whitney regular
stratification of $(\mfk g^*\oplus\mfk g)^c$ pulls back to a Whitney regular
stratification of $(P\times\mfk g)^c$. The strata of $(P\times\mfk g)^c$
are the nonempty submanifolds
\begin{equation*}\begin{split}
(P\times\mfk g)_{(K)}&=\widetilde{J}^{-1}\bigl((\mfk g^*\oplus
\mfk g)_{(K)}\bigr)\\
&=\set{(p,\xi)\in P\times\mfk g}{\mbox{$G_{J(p),\xi}$ is conjugate to $K$}}
\end{split}\end{equation*}
where $(K)\in\mcl I^c$.

Let $\mcl T$ denote the $G$-invariant subbundle of the vector
bundle $TP$ with fiber over $p\in P$ equal to the tangent space at $p$
to the orbit of $G$ through $p$, that is $\mcl T_p=\mfk g.p$. The
subbundle $\mcl T$ is the image of the map
$I:P\times\mfk g\rightarrow TP$ defined by $(p,\xi)\mapsto\xi.p$.
Since $G$ is acting freely on $P$, the map $I$ is a
$G$-equivariant vector bundle isomorphism between $P\times\mfk g$
(considered as a bundle over $P$) and $\mcl T$.

Let $\mcl K$ denote the $G$-invariant subbundle of the vector bundle
$TP$ with fiber over $p\in P$ given by $\mcl K_p =\ker dJ(p)$. Note
that the Hamiltonian vector field of any $G$-invariant Hamiltonian
on $P$ takes values in $\mcl K$. Let $\mcl T^c =\mcl T\cap\mcl K$, a
$G$-invariant subset of $\mcl K$. For each $(K)\in\mcl I^c$ define the
set $\mcl T^c_{(K)}$ by
\begin{equation*}
\mcl T^c_{(K)}=\bset{\xi.p\in\mcl T^c}{\mbox{$\bigl(J(p),
\xi\bigr)$ has type $(K)$}}.
\end{equation*}

Let $\mcl T^o$ and $\mcl K^o$ denote the $G$-invariant vector
subbundles of $T^*P$ with fibers
\begin{equation*}
\mcl T_{p}^{o}=\onm{ann}(\mcl T_p),\quad\mcl K_{p}^{o}=
\onm{ann}(\mcl K_p).
\end{equation*}
The subbundle $\mcl K^o$ is the image of the map
$I^o:P\times\mfk g\rightarrow T^*P$ defined by $(p,\xi)\mapsto
dJ_\xi(p)$. Moreover, since $G$ is acting freely on $P$, the map
$I^o$ is a $G$-equivariant vector bundle isomorphism between $P\times
\mfk g$ (considered as a bundle over $P$) and $\mcl K^o$. Note that
any $G$-invariant Hamiltonian $P$ is such that $dH$ has values in
$\mcl T^o$. Let $\mcl K^{oc} =\mcl T^o\cap\mcl K^o$, a
$G$-invariant subset of $\mcl T^o$. For each $(K)\in\mcl I^c$ define
the set $\mcl K^{oc}_{(K)}$ by
\begin{equation*}
\mcl K^{oc}_{(K)}=\bset{dJ_\xi(p)\in\mcl K^{oc}}{
\mbox{$\bigl(J(p),\xi\bigr)$ has type $(K)$}}.
\end{equation*}

\begin{proposition}\mbox{}
\begin{enumerate}
\item 
The set $\set{\mcl T^c_{(K)}}{(K)\in\mcl I^c}$
is a Whitney regular
stratification of $\mcl T^c$ and the diffeomorphism
$I:P\times\mfk g\rightarrow\mcl T$ restricts to a stratum preserving
bijection of $(P\times\mfk g)^c$ to $\mcl T^c$.
\item
The set $\set{\mcl K^{oc}_{(K)}}{(K)\in\mcl I^c}$
is a Whitney regular stratification of $K^{oc}$ and the diffeomorphism
$I^o:(P\times\mfk g)\rightarrow\mcl K^o$ restricts to a stratum preserving
bijection of $(P\times\mfk g)^c$ to $\mcl K^{oc}$.
\end{enumerate}
\end{proposition}

\begin{proof}
For part~1, equivariance of $J$ implies that
$dJ(p)(\xi.p)=-\onm{coad}_\xi J(p)$ and so
\begin{equation*}
I(p,\xi)\in\mcl T^c\ \Leftrightarrow \ \xi.p\in\mcl K\ \Leftrightarrow\ 
0=-\onm{coad}_\xi J(p)\ \Leftrightarrow\ (p,\xi)\in (P\times\mfk g)^c.
\end{equation*}
Hence $I$ is a diffeomorphism which maps $(P\times\mfk g)^c$
bijectively to $\mcl T^c$. The strata $\mcl T^{c}_{(K)}$ are the
images under $I$ of the strata of $(P\times\mfk g)^c$ and so define a
Whitney regular stratification of $T^c$. The map $I$ is stratum
preserving by construction. Part~2 is similar. \qed\end{proof}

Let $H:P\rightarrow\mbb R$ be a $G$-invariant Hamiltonian function on
$P$, $dH:P\rightarrow\mcl T^o\subset T^*P$ its derivative and
$X_{H}:P\rightarrow\mcl K\subset TP$ the corresponding
$G$-equivariant Hamiltonian vector field.
Denote by $\psi:P\times\mfk g\rightarrow TP$ the map $(p,\xi)\mapsto
X_H(p)-\xi.p$, and let $\psi_{(K)}$ denote the restriction of $\psi$
to the submanifold $(P\times\mfk g)^{c}_{(K)}$ of $P\times\mfk g$.
Denote by $\psi^o:P\times\mfk g\rightarrow T^*P$ 
the map $(p,\xi)\mapsto dH(p)-dJ_\xi(p)$
and let $\psi^{o}_{(K)}$ denote its restriction to
$(P\times\mfk g)^{c}_{(K)}$. Note that $\psi_{(K)}$ maps into $\mcl K$
and $\psi^{o}_{(K)}$ into $\mcl T^o$. The following lemma is a
direct consequence of the definitions.

\begin{lemma}
The vector field $X_H$ has a relative equilibrium at $p_e$ if and only
if the following equivalent conditions hold:
\begin{enumerate}
\item $X_H(p_e)\in\mcl T^c$;
\item $dH(p_e)\in\mcl K^{oc}$;
\item There exists $\xi_e\in\mfk g$ such that
$(p_e,\xi_e)\in (P\times\mfk g)^c$ and $\psi(p_e,\xi_e)=0$;
\item There exists $\xi_e\in\mfk g$ such that
$(p_e,\xi_e)\in (P\times\mfk g)^c$ and $\psi^o(p_e,\xi_e)=0$.
\end{enumerate}
In statements 3 and 4 the element $\xi_e$ is the generator of the relative
equilibrium. Moreover, if $p_e$ is a relative equilibrium of $X_H$ with
generator $\xi_e$, then
\begin{multline*}
\mbox{$p_e$ has type $(K)$}
\ \Leftrightarrow\ (p_e,\xi_e)\in (P\times\mfk g)^c_{(K)}\\
\Leftrightarrow\ X_H(p_e)\in\mcl T^c_{(K)}
\ \Leftrightarrow\ dH(p_e)\in\mcl K^{oc}_{(K)}.
\end{multline*}
\end{lemma}

The following lemma and the accompanying definition identify
what we mean by a transversal relative equilibrium. As we shall show,
transversality is generic in the class of symmetric Hamiltonian
systems.

\begin{lemma}\lb{transcdns} 
Let $p_e$ be a relative equilibrium of $X_H$ with generator $\xi_e$ and 
momentum $\mu_e=J(p_e)$, and
let $K_e=G_{\mu_e}\cap G_{\xi_e}$. Then
the following are equivalent:
\begin{enumerate}
\item 
$X_{H}:P\rightarrow\mcl K$ is transversal
to $\mcl T^c_{(K_e)}$ at $p_e$;
\item 
$dH:P\rightarrow\mcl T^o$ is transversal
to $\mcl K^{oc}_{(K_e)}$ at $p_e$;
\item 
$\psi_{(K_e)}:(P\times\mfk g)^c_{(K_e)}\rightarrow\mcl K$ is transversal to
the zero section of $\mcl K$ at~$(p_e,\xi_e)$;
\item 
$\psi^o_{(K_e)}:(P\times\mfk g)^c_{(K_e)}\rightarrow\mcl T^o$ is transversal to
the zero section of $\mcl T^o$ at~$(p_e,\xi_e)$.
\end{enumerate}
\end{lemma}

\begin{proof}
The symplectic form $\omega$ on $P$ defines a vector bundle isomorphism
$\omega^\flat:TP\rightarrow T^*P$ which converts
part~1 into part~2, since
\begin{equation*}
\omega^\flat\circ X_H=dH,\quad\omega^\flat(\mcl K_e)=T^o,\quad
\omega^\flat(\mcl T^c_{(K_e)})=\mcl K^{oc}_{(K_e)}.
\end{equation*}
Thus~(1) is equivalent to~(2) and similarly~(3) is equivalent
to~(4), so we need only show the equivalence of~(1) and ~(3).

Generally, as is easily verified, given a section $X$ of a vector
bundle $\pi:E\rightarrow P$ and a mapping $f$ from a manifold
$M$ to $E$ over $f_0:M\rightarrow P$, the condition that $X$ and $f$
are transversal is equivalent to the condition that $f-X\circ f_0$ is
transversal to the zero section $Z(E)$ of $E$. When
applied to the diagram
\begin{equation*}\begin{diagram}
\node{}\node{\mcl K}\\
\node{(P\times\mfk g)^c_{(K_e)}}\arrow{e,b}{(p,\xi)\mapsto p}\arrow{ne,t}{I}
\node{P}\arrow{n,r}{X_H}
\end{diagram}\end{equation*}
this shows that~(3) is equivalent to the statement that
$I|(P\times\mfk g)^c_{(K_e)} =\mcl T^c_{(K_e)}$ is transversal to $X_H$,
and so is equivalent to~(1).
\qed\end{proof}

\begin{definition}\lb{defn}
A relative equilibrium $p_e$ is said to be {\bfseries transversal} if
the equivalent conditions in Lemma~\ref{transcdns} hold.
\end{definition}

Let $\mcl E\subset P$ denote the set of all relative equilibria of $X_H$
and $\mcl E_{(K)}$ the subset consisting of those of type $K$.
%We also have
%the subsets $\mcl E/G$ and $\mcl E_{(K)}/G$ of $P/G$.

\begin{theorem}\lb{mainth1}
Let $p_e$ be a relative equilibrium of $X_H$ with
generator $\xi_e$ and momentum $\mu_e=J(p_e)$, and let
$K_e=G_{\mu_e}\cap G_{\xi_e}$.
If $p_e$ is transversal then there exists a $G$-invariant
open neighborhood $U$ of $p_e$ in $P$ such that
\begin{enumerate}
\item 
every relative equilibrium of $X_H$ in $\mcl E\cap U$ is transversal;
\item
for each $(K)\in\mcl I^c(K_e)$ the subset $\mcl E_{(K)}\cap U$ is a
submanifold of dimension $\dim G+2\dim Z(K)-\dim K$ and its quotient
by the action of $G$, $(\mcl E_{(K)}\cap U)/G$ is a manifold of
dimension $2\dim Z(K)-\dim K$;
\item 
the sets $\set{\mcl E_{(K)}\cap U}{(K)\in\mcl I^c(K_e)}$
and $\set{(\mcl E_{(K)}\cap U)/G}{(K)\in\mcl I^c(K_e)}$
are Whitney regular stratifications of $\mcl E\cap U$ and
$(\mcl E\cap U)/G$, respectively.
\end{enumerate}
\end{theorem}

\begin{proof}
Since $\mcl E=X_H^{-1}(\mcl T^c)$ the statements about $\mcl E$ follow
directly from the transversality of $X_H$ at
$p_e$ to the Whitney regular stratification of $\mcl T^c$ described above.
The results on $\mcl E/G$ follow from those on $\mcl E$ and the fact that
$G$ acts freely on $P$.
\qed\end{proof}

The possible types of a transversal relative equilibrium are severely
restricted by the following corollary. For example, if $G$ has no
center then no relative equilibria can have type $(G)$.

\begin{corollary}
If $p_e$ is a transversal relative equilibrium then
\begin{equation*}\dim Z(K_e)\ge\frac12\dim K_e.\end{equation*}
\end{corollary}

\begin{proof}
This is immediate since
$\onm{dim}\mcl E_{(K_e)}/G =2\dim Z(K_e)-\dim K_e$
is non-negative.
\qed\end{proof}

We next give a more detailed description of the local structure of the
set $\mcl E$. First we recall some properties of Whitney stratifications,
as described in \ct{GoreskyMMacphersonR-1988.1}. Let $M$ be a manifold,
$Z\subset M$ a Whitney stratified space, $p\in Z$, and $S$ the stratum
through $p$. Let $N^\prime$ be a submanifold of $M$ transverse to $S$
at $p$, meaning that $T_pM=T_pS\oplus T_pN^\prime$.
Choose a Riemannian metric on $N^\prime$. Then for
small enough $\delta$ the {\em normal slice \/} $N= N^\prime\cap
Z\cap B_\delta(p)$ and the {\em link} $L= N^\prime\cap Z\cap
\partial B_\delta(p)$ have topological types independent of $\delta$,
the Riemannian metric and $N^\prime$. Moreover $N$ is a cone over $L$
and locally $Z$ is the product of $N$ and $S$. These result follow from
the Thom isotopy lemma.

Suppose that near $p$ the set $Z$ is a smooth stratum preserving
embedding of the product of a fixed stratified cone $C$ in some vector space
$\mbb E$ and an open subset $U$ of Euclidean space. Proposition~\ref{prop1}
shows that this applies to points in $(\mfk g^*\oplus\mfk g)^c$.
Without loss of generality we may assume that
the cone linearly spans $\mbb E$. Then we will say that $Z$ 
has (smooth) singularity-type $C$ at $p$. We make the following elementary
observations, all based on more or less transparent applications of
the implicit function theorem:
\begin{enumerate}
\item 
Given any transverse submanifold $N^\prime$, the normal space $N$ is
the image of a smooth embedding of $C$ into $N^\prime$. Conversely, if
$N$ is such for some such $N^\prime$ then $Z$ has singularity type $C$
at $p$.
\item 
Suppose $M^\prime$ is a submanifold of $M$ and $Z\subseteq M^\prime$.
Take any transverse submanifold $N^{\prime\prime}$ in $M^\prime$ to
the stratum $S$ and extend it to a transverse submanifold $N^\prime$
to $S$ within $M$. By Item~(1) just above the normal space $N$ is the
image of $C$ by a smooth embedding $\iota :\mbb E\rightarrow N^\prime$,
say. Then $\iota(\mbb E)$ is a submanifold with tangent space contained in
$N^{\prime\prime}$, and consequently $\iota$ may be deformed so that
its image is contained in $N^{\prime\prime}$, all the while fixing its
values on $C$. The point of this is that $Z$ also has singularity-type $C$
at $p$ when it is regarded as a subset of $M^\prime$.
\item 
If $N$ is some manifold and $f:N\rightarrow M$ is a smooth map such that
$f(n)=p$ and $f$ is transversal to the stratum $S$ at $n$, then
$f^{-1}(Z)$ is a stratified space with the same singularity type at
$n$ as that of $Z$ at $p$. For this we invoke smooth coordinates at
$z$ such that $Z$ locally becomes $C\times S\times\{0\}\subset\mbb
E\times S\times\mbb F$, where $\mbb F$ is some vector
space. Then $f$ is transversal to $\mbb E\times S\times\{0\}$, so
we may replace $f$ by its restriction to $f^{-1}(\mbb E\times
S\times\{0\})$, thereby removing $\mbb F$. Let $\pi_1:\mbb
E\times S\rightarrow\mbb E$ be the projection. Again by
transversality, $\pi_1\circ f$ is a submersion at $n$ and so, by local
diffeomorphism of $N$ near $n$, $\pi_1\circ f$ becomes the projection
$(x,y)\mapsto x$, and consequently $f^{-1}(Z)$ also has
singularity-type $C$ at $n$.
\end{enumerate}
Our point is just that these singularity theory results are obtainable
in the smooth category without the use of the Thom isotopy lemma
provided the smooth structure of the stratified spaces is
apriori known, as in Proposition~\ref{prop1}.

\begin{theorem}\lb{mainth2}
The stratified spaces $\mcl E\cap U$ and $(\mcl E\cap U)/G$ of
Theorem~\ref{mainth1} both have singularity type
$\bigl(\mfk l_e^*\oplus\mfk l_e\bigr)^c$ at $p_e$ and $\bar p_e = G.p_e$,
respectively, where $\mfk l_e$ is the Lie algebra of $K_e/Z(K_e)$.
\end{theorem}

\begin{proof}
By Propositions~\ref{prop1} and~\ref{prop2}, the set $(\mfk g\oplus\mfk g)^c$
has singularity type
$\bigl(\mfk l_e^*\oplus\mfk l_e\bigr)^c$ at $(J(p_e),\xi_e)$.
Since $J\times\Id:P\times\mfk g\rightarrow\mfk
g^*\times\mfk g$ is a submersion, the singularity type of
$(P\times\mfk g)^c=(J\times\Id)^{-1}(\mfk g^*\oplus\mfk g)^c$ at
$(p_e,\xi_e)$ is also $\bigl(\mfk l_e^*\oplus\mfk l_e\bigr)^c$.
The set $T^c$ also has this singularity type at the point $(p_e,\xi_e.p_e)$,
since it is mapped diffeomorphically to $(P\times\mfk g)^c$ by $I$.
Regarding $T^c$ as a subset of $\mcl K$ leaves its singularity-type
unchanged. By transversality, and since $\mcl
E=X_H^{-1}(\mcl T^c)$, this is also the singularity type of $\mcl E$
at $p_e$. Since $G$ acts freely on $P$ the set $\mcl E$ is locally
isomorphic to the product of $\mcl E/G$ and $G$ and so $\mcl E/G$
has the same singularity type at $\bar p_e$ as $\mcl E$ does at $p_e$. 
\qed\end{proof}

We end this section with a result that shows that all the relative equilibria
of `most' $G$-invariant Hamiltonians are transversal.
Recall that $G$ is a connected, compact Lie group acting freely and
symplectically on the symplectic manifold $(P,\omega)$. Let $C_{G}^\infty(P)$
denote the set of all $G$-invariant $C^\infty$ functions on $P$ and
$C^\infty(P/G)$ the set of all $C^\infty$ functions on the orbit
space $P/G$. Pulling back functions by the orbit map $\pi : P\rightarrow
P/G$ defines a bijection
$\pi^* : C^\infty(P/G)\cong C_{G}^\infty(P/G)$.
We give $C^\infty(P/G)$ the Whitney $C^\infty$ topology and $C_{G}^\infty(P)$ the
isomorphic topology induced by $\pi^*$.

\begin{theorem}\lb{mainth3}
There exists a dense subset $\mcl H\subset C_{G}^\infty(P)$ such
that if $H\in\mcl H$ then every relative equilibrium of $H$ is transversal.
\end{theorem}

\begin{proof}
For any $G$-invariant function $H$ on $P$ let $\bar H$ denote the
quotient function on $P/G$. The derivative $dH:P\rightarrow\mcl
T^o\subset T^*P$ is $G$-equivariant and so descends to a mapping from
$P/G$ to $\mcl T^o/G$ which we will denote by
$(dH)\bar{\phantom{\rule{3pt}{5pt}}}$. Since $G$ is acting freely on
$P$ the orbit space $P/G$ is a smooth manifold and $\mcl T^o/G$ is a
smooth vector bundle over $P/G$ which can be identified with the
cotangent bundle $T^*(P/G)$. Under this identification
$(dH)\bar{\phantom{\rule{3pt}{5pt}}}$ becomes equal to $d\bar H$.

The Hamiltonian $H$ has a relative equilibrium at $p_e\in P$ if and only
if $dH(p_e)\in\mcl K^{oc}$ and so if and only if
$d\bar H(\bar p_e)\in\mcl K^{oc}/G$. The set
$\set{\mcl K^{oc}_{(K)}/G}{(K)\in\mcl I^c}$ is a Whitney regular
stratification of $\mcl K^{oc}/G\subset T^*(P/G)$ and $p_e$ is a transversal
relative equilibrium if and only if $d\bar H : P/G\rightarrow T^*(P/G)$
is transversal to this stratification at $\bar p_e$. It now follows from the
Thom transversality theorem that the set of functions $\bar H$ on $P/G$
for which the derivatives $d\bar H$ are everywhere transversal to this
stratification is dense in $C^\infty(P/G)$. Pulling this back to
$C_{G}^{\infty}(P)$ gives the required set $\mcl H$.
\qed\end{proof}

Thus for $H$ in the set $\mcl H$ given by the Theorem the set of all relative
equilibria $\mcl E\subset P$ is Whitney stratified by the types of the
relative equilibria, the manifold $\mcl E_{(K)}$ of relative equilibria of
type $(K)$ has dimension $\dim G+2\dim Z(K)-\dim K$, and the
singularity-type of $\mcl E$ at a point in $\mcl E_{(K)}$ is
$(\mfk l^*\oplus\mfk l)^c$ where $\mfk l$ is the Lie algebra of
$K/Z(K)$.

\section{Linearization of a Transversal Relative Equilibrium}\lb{linear}

In this section we give an alternative form of the
transversality condition using a normal form for the linearization
of a vector field at a relative equilibrium due to~\ct{PatrickGW-1998.1}.
This is then used to describe the tangent spaces of the strata of $\mcl E/G$
and to prove a generalization and partial converse of a result of
\ct{PatrickGW-1998.1} on the symplectic properties of these strata.

Recall that $\mcl T_{p_e}=\mfk g.p_e$ and $\mcl
K_{p_e}=\ker dJ(p_e)$, and decompose $T_{p_e}P$ as a direct sum
\begin{equation}\lb{decomp}
T_{p_e}P=T_0\oplus N_1\oplus N_0\oplus T_1
\end{equation} 
as follows. Let $T_0=\mfk g_{\mu_e}.p_e$.
Let $\mfk g_{\mu_e}^\perp$ be a $G$-invariant complement to
$\mfk g_{\mu_e}$ in $\mfk g$ and set $T_1=\mfk g_{\mu_e}^\perp.p_e$.
Choose a complement $N_1$ to $T_0$ within $\mcl K_{p_e}$ and a complement
$N_0$ to $\mcl T_{p_e}\oplus\mcl K_{p_e}$ within $T_{p_e}P$.
The subspace $N_1$ is a symplectic
normal space at $p_e$ and can be identified with the tangent space to
the reduced phase space $P_{\mu_e}$ at $\bar p_e= G.p_e$.
The subspace $N_0$ can be identified with $\mfk g_{\mu_e}^*$, regarded as
a subspace of $\mfk g^*$, via the map $dJ(p_e)|N_0$. Since $p_e$ is regular,
we have $T_0\cong\mfk g_{\mu_e}$ and $T_1\cong\mfk g_{\mu_e}^\perp$, and
so~\rf{decomp} becomes
\begin{equation}\lb{decomp2}
T_{p_e}P\cong\mfk g_{\mu_e}\oplus T_{\bar p_e}P_{\mu_e}\oplus
\mfk g_{\mu_e}^*\oplus\mfk g_{\mu_e}^\perp.
\end{equation}
The symplectic form with respect to
this decomposition is the product of the reduced form on
$T_{\bar p_e}P_{\mu_e}$, the canonical form on $\mfk
g_{\mu_e}\oplus\mfk g_{\mu_e}^*$, and the Kostant-Souriau form on
$\mfk g_{\mu_e}^\perp=T_{\mu_e}(G.\mu_e)$. The map
$dJ(p_e)$ has the explicit form
\begin{equation}\lb{dJform}
dJ(p_e)(\xi_0\oplus w\oplus\mu_0\oplus\xi_1)=\mu_0+\onm{coad}_{\xi_1}\mu_e.
\end{equation}

A relative equilibrium $p_e$ of $H$ is an
equilibrium point of $X_{H_{\xi_e}}$ where $H_{\xi_e}=H-J_{\xi_e}$.
\ct{PatrickGW-1998.1} shows that the complements $N_0$ and $N_1$ can be
chosen so that, with respect to the decomposition~\rf{decomp2}, the
linearization of $X_{H_{\xi_e}}$ at $p_e$ has the form
\begin{equation}\lb{lin}
dX_{H_{\xi_e}}(p_e)= 
\left[\begin{array}{cccc}
-\onm{ad}_{\xi_e}\!\!|\mfk g_{\mu_e} &C^*&D&0\\
0&dX_{H_{\mu_e}}(\bar p_e)&C&0\\
0&0&-\onm{coad}_{\xi_e}\!\!|\mfk g_{\mu_e}^* &0\\
0&0 &0&-\onm{ad}_{\xi_e}\!\!|\mfk g_{\mu_e}^\perp
\end{array}\right]\end{equation}
where $ H_{\mu_e}$ is the induced Hamiltonian on the reduced phase
space $P_{\mu_e}$, $X_{H_{\mu_e}}$ is the associated Hamiltonian
vector field and $dX_{H_{\mu_e}}(\bar p_e)$ is its linearization at
the equilibrium point $\bar p_e$. The operator $D:\mfk
g_{\mu_e}^{*}\cong N_0\rightarrow T_0\cong \mfk g_{\mu_e}$ describes
the {\it drift\/} along the group orbit and $C:\mfk g_{\mu_e}^{*}\cong
N_0\rightarrow N_1$ describes interactions between the reduced
dynamics and the motion along group orbits due to possible 1-1
resonances between these two motions. If the spectrum of
$\onm{ad}_{\xi_e}$ is distinct from the spectrum of
$dX_{H_{\mu_e}}(\bar p_e)$ then such interactions do not occur to
first order and $N_0$ and $N_1$ can be chosen so that $C=0$. In
general the dual operator $C^*$ is regarded as a map from $N_1$
(instead of $N_1^*$) by using the symplectic form to identify $N_1$
with its dual. The operators $C$, $C^*$ and $D$ enjoy the following
properties:
\begin{gather}
dX_{H_{\mu_e}}(\bar p_e)C=-C\onm{coad}_{\xi_e},\lb{cprop1}\\
C^*dX_{H_{\mu_e}}(\bar p_e)=-\onm{ad}_{\xi_e}C^*,\lb{cprop2}\\
\onm{ad}_{\xi_e}D=D\onm{coad}_{\xi_e},\lb{cprop3}
\end{gather}
and the operator $D$ is symmetric.

As in the theorems of the previous section we will set
$K_e=G_{\mu_e}\cap G_{\xi_e}$, denote by $\mfk k_e$ the Lie algebra of
$K_e$, $\mfk z_e$ the center of $\mfk k_e$, and $\mfk l_e=\mfk
k_e/\mfk z_e$. As a consequence of~\rf{cprop1}, $C$ maps $\mfk k_e^*$
into $\onm{ker}dX_{H_{\mu_e}}(\bar p_e)$, since $\onm{coad}_{\xi_e}$
is zero on $\mfk k_e^*$. Equation~\rf{cprop2} similarly implies that
$C^*$ maps $\onm{ker}dX_{H_{\mu_e}}(\bar p_e)$ into $\mfk k_e$ while
equation~\rf{cprop3} implies that $D$ maps $\mfk k_e^*$ into $\mfk k_e$.
The next theorem characterizes transversal relative equilibria
in terms of the above normal form for $dX_{H_{\xi_e}}(p_e)$.
Recall that $p_e$ is said to be {\it nondegenerate} if
$dX_{H_{\mu_e}}(\bar p_e)$ is invertible.
\begin{theorem}\lb{normaltransv}
The relative equilibrium $p_e$ is transversal if and only if all the
following three conditions hold:
\begin{enumerate}
\item either $p_e$ is nondegenerate or $0$ is a
semisimple eigenvalue of $dX_{H_{\mu_e}}(\bar p_e)$;
\item
$C$ maps $\mfk z_e^*\subset\mfk k_e^*\subset\mfk g_{\mu_e}^*$ onto
$\ker dX_{H_{\mu_e}}(\bar p_e)$;
\item
$C^*\bigl(\ker dX_{H_{\mu_e}}(\bar p_e)\bigr)+D\bigl(\ker C\cap 
\mfk z_e^*\bigr)+\mfk z_e=\mfk k_e$.
\end{enumerate}
\end{theorem}

\begin{proof}

From~\ct{PatrickGW-1995.1}, page~409, there is the following formula
for the derivative of $\psi:P\times\mfk g\rightarrow TP$ at $(p_e,\xi_e)$
\begin{equation}\lb{Tpsi}
T_{(p_e,\xi_e)}\psi(v,\eta)=\bigl(v,dX_{H_{\xi_e}}(p_e)v-\eta.p_e\bigr),
\end{equation}
where on the right side $T_{0_{p_e}}(TP)$ has been decomposed
into horizontal and vertical parts.
By Lemma~\ref{transcdns}, $p_e$ is transversal if and only if 
$\psi_{(K_e)}:(P\times\mfk g)^c_{(K_e)}\rightarrow\mcl K$
is transversal to
the zero section of $\mcl K$ at $(p_e,\xi_e)$, which is equivalent to
\begin{equation*}
\ker dJ(p_e)\subseteq\bset{dX_{H_{\xi_e}}(p_e)v-\eta.p_e}{(v,\eta)\in 
T_{(p_e,\xi_e)}(P\times\mfk g)_{(K_e)}^c},
\end{equation*}
while from Propositions~\ref{prop1} and~\ref{prop2},
\begin{multline*}
T_{(p_e,\xi_e)}(P\times \mfk g)_{(K_e)}^c\\
=\bset{(v,\eta)\in T_{p_e}P\times\mfk g}
{\bigl(dJ(p_e)v,\eta\bigr)\in T_{(\mu_e,\xi_e)}(\mfk g^*\oplus
\mfk g)^c_{(K_e)}}\\
=\bset{(v,\eta)\in T_{p_e}P\times\mfk g}
{\bigl(dJ(p_e)v,\eta\bigr)\in
(\mfk z_e^* \oplus \mfk z_e)\oplus\mfk g.(\mu_e,\xi_e)}.
\end{multline*}
Using~\rf{dJform} and~\rf{lin} transversality is equivalent to
\begin{multline}
\mfk g_{\mu_e}\oplus N_1=\bigl\{\,\bigl(
-\onm{ad}_{\xi_e}\xi_0+C^*w+D\mu_0-\eta_0,
 dX_{H_{\mu_e}}(\bar p_e)w+C\mu_0\bigr):\\
\xi_0,\ \eta_0\in\mfk g_{\mu_e},\ w\in T_{\bar p_e}P_{\mu_e},\ 
 \mu_0\in\mfk g_{\mu_e}^{*}
\mbox{ and conditions~\rf{cond1}--\rf{cond4} hold}\,\bigr\},\lb{transv0}
\end{multline}
the conditions~\rf{cond1}--\rf{cond4} being that there exist
$\xi_1, \eta_1 \in \mfk g_{\mu_e}^{\perp},\zeta \in \mfk z_e^*,z\in \mfk z_e$
and $\tilde\xi\in\mfk g$ such that
\begin{gather}
\onm{coad}_{\xi_e}\mu_0=0,\lb{cond1}\\
\eta_1=-\onm{ad}_{\xi_e}\xi_1,\lb{cond2}\\
\mu_0+\onm{coad}_{\xi_1}\mu_e=\zeta+\onm{coad}_{\tilde\xi}\mu_e,\lb{cond3}\\
\eta_0+\eta_1=z+\onm{ad}_{\tilde\xi}\xi_e.\lb{cond4}
\end{gather}

Condition~\rf{cond1} gives
$\mu_0\in\onm{ann}\bigl(\mfk g_{\mu_e}.\xi_e\bigr) = \mfk k_e^*$.
Inserting~\rf{cond2} into~\rf{cond4}, one
sees that~\rf{cond1}--\rf{cond4} are equivalent to the conditions that
$\mu_0\in\mfk k_e^*$ and there exist $\zeta\in\mfk z_e^*$, $z\in \mfk z_e$ and
$\tilde\xi\in\mfk g$ such that
\begin{gather}
\mu_0=\zeta+\onm{coad}_{\tilde\xi-\xi_1}\mu_e,\lb{cond5}\\
\eta_0=z+\onm{ad}_{\tilde\xi-\xi_1}\xi_e.\lb{cond6}
\end{gather}
In~\rf{cond5} $\mu_0$ and $\zeta$ lie in $\mfk k_e^*$ which is
orthogonal to the image of $\onm{coad}$ and so~\rf{cond5} is equivalent
to $\mu_0\in \mfk z_e^*$ and $\tilde\xi-\xi_1\in\mfk
g_{\mu_e}$. Since the right side of~\rf{cond6} is then clearly in
$\mfk g_{\mu_e}$, after substitution of~\rf{cond6} into~\rf{transv0},
transversality becomes equivalent to
\begin{multline}
\mfk g_{\mu_e}\oplus N_1
=
\bigl\{\,\bigl(C^*w+D\mu_0+g_{\mu_e}.\xi_e + \mfk z_e,
 dX_{H_{\mu_e}}(\bar p_e)w+C\mu_0\bigr):\\
w \in T_{\bar p_e}P_{\mu_e},\ \mu_0 \in \mfk z_e^*\,\bigr\}.
 \lb{transv1} 
\end{multline}

Let $E_0$ be the generalized eigenspace of $dX_{H_{\mu_e}}(\bar p_e)$
associated to the eigenvalue $0$ and let $E_1$ be the sum of the other
generalized eigenspaces, so that $N_1=E_0\oplus E_1$. This
decomposition is preserved by $dX_{H_{\mu_e}}(\bar p_e)$ and for $j =
0,1$ we denote the restriction of $dX_{H_{\mu_e}}(\bar p_e)$ to $E_j$
by $L_j$. It follows from \rf{transv1} that for $p_e$ to be
transversal we must have $\onm{image}L_0+\onm{image}C\supseteq
E_0$. Since $\onm{image}C \subset \onm{ker}dX_{H_{\mu_e}}(\bar p_e) =
\onm{ker}L_0 \subset E_0$ this implies $\onm{image}L_0 + \onm{ker}L_0
= E_0$. For a nilpotent operator, such as $L_0$, an application of
the Jordan normal form shows that this is only possible if $L_0 = 0$.
It follows that if $p_e$ is transversal then either it is
nondegenerate or $0$ is a semisimple eigenvalue of
$dX_{H_{\mu_e}}(\bar p_e)$.

So in the proving of either direction of the theorem we may assume
$L_0=0$. But then the transversality condition~\rf{transv1} reduces to
\begin{multline}\lb{300}
\mfk k_e\oplus\onm{ker}dX_{H_{\mu_e}}(\bar p_e)=\bigl\{\,\bigl(C^*w+D\mu_0
 +\mfk g_{\mu_e}.\xi_e+\mfk z_e, C\mu_0\bigr): \\
 w\in\ker dX_{H_{\mu_e}}(\bar p_e),\ \mu_0\in\mfk z_e^*\,\bigr\}.
\end{multline}
As $C^*$ maps $\onm{ker}dX_{H_{\mu_e}}(\bar p_e)$ into $\mfk k_e$,
while $D$ maps $\mfk z_e^*$ into $\mfk k_e$, $\mfk z_e\subseteq\mfk k_e$,
and $\mfk g_{\mu_e}.\xi_e=\mfk k_e^\perp\cap\mfk g_{\mu_e}$,~\rf{300} is
equivalent to the pair of statements 2 and 3 in the Theorem. \qed
\end{proof}

Let $\bar D:\mfk z_e^* \rightarrow \mfk l_e \cong \mfk k_e/\mfk z_e$ denote
the mapping obtained by the restricting the drift operator $D$ to
$\mfk z_e^*$ and then composing with the projection from $\mfk k_e$ to
$\mfk l_e\cong\mfk k_e/\mfk z_e$.
\begin{corollary}\mbox{}
\begin{enumerate}
\item If $p_e$ is nondegenerate then $p_e$ is transversal if and only
if $\bar D$ is surjective.
\item If $p_e$ is transversal then
\begin{equation*}
\dim\ker dX_{H_{\mu_e}}(\bar p_e)\leq 2\dim Z(K_e)-\dim K_e.\end{equation*}
\end{enumerate}
\end{corollary}

\begin{proof}
If $p_e$ is nondegenerate then $\ker dX_{H_{\mu_e}}(\bar p_e) = \{0\}$,
$\ker C = \mfk g_{\mu_e}^{*}$ and the transversality conditions reduce
to $D(\mfk z_e^*)+\mfk z_e=\mfk k_e$, as required.

The second statement follows immediately from the third transversality
condition in Theorem~\ref{normaltransv}.\qed
\end{proof}

We note that if $K_e$ is a maximal torus then $\mfk z_e = \mfk k_e$ and so 
$\bar D$ is trivially surjective and nondegeneracy implies
transversality. However in general a non-trivial nondegeneracy condition
must be satisfied by the drift operator for a relative equilibrium to
be transversal.
 
The normal form~(\ref{lin}) also enables us to give a description of the
tangent spaces to the manifolds $\mcl E_{(K_e)}$ at transversal relative
equilibria.

\begin{theorem}\lb{tanspace}
If the relative equilibrium $p_e$ is transversal then, with respect to the
decomposition~\rf{decomp2},
\begin{multline*}
T_{p_e}\mcl E_{(K_e)}=\bigl\{\,
\xi_0\oplus w\oplus \mu_0\oplus\xi_1:\xi_0\in\mfk g_{\mu_e},\ 
 \xi_1\in\mfk g_{\mu_e}^{\perp},\ \mu_0\in\mfk z_e^*\cap\ker C,\\
w \in \ker dX_{H_{\mu_e}}(\bar p_e),\ 
C^*w+D\mu_0\in\mfk z_e\,\bigr\}.
\end{multline*}\end{theorem}

\begin{proof}
Since $\mcl E_{(K_e)}$ is the projection to $P$ of
$\psi_{(K_e)}^{-1}\bigl(Z(TP)\bigr)\subset( P \times \mfk g)^{c}_{(K_e)}$,
and in view of~\rf{Tpsi}, we have
\begin{multline*}
T_{p_e}\mcl E_{(K_e)}=\bigl\{v\in T_{p_e}P: 
 dX_{H_{\xi_e}}(p_e)v-\eta.p_e=0 \\
\mbox{for some $\eta\in\mfk g$ such that
$(v,\eta) \in T_{(p_e,\xi_e)}(P \times \mfk g)^{c}_{(K_e)}$}\bigr\},
\end{multline*}
which, in terms of the local normal form~\rf{lin}, becomes
\begin{equation*}
T_{p_e}\mcl E_{(K_e)}=\set{\xi_0\oplus w\oplus\mu_0\oplus\xi_1}{
\mbox{the conditions directly below hold}},
\end{equation*}
the conditions being that there exists $\eta_0\in\mfk g_{\mu_e}$ and
$\eta_1\in\mfk g_{\mu_e}^\perp$ such that
\begin{gather*}
-\onm{ad}_{\xi_e}\xi_0+C^*w+D\mu_0-\eta_0=0,\quad
dX_{H_{\mu_e}}(\bar p_e)w+C\mu_0=0,\\
\onm{coad}_{\xi_e}\mu_0=0,\quad
-\onm{ad}_{\xi_e}\xi_1-\eta_1=0,
\end{gather*}
subject to the constraints that there exist
$\zeta\in \mfk z_e^*, z\in \mfk k_e$ and $\tilde\xi\in\mfk g$ such that
\begin{gather*}
\mu_0+\onm{coad}_{\xi_1}\mu_e=\zeta+\onm{coad}_{\tilde\xi}\mu_e,\quad
\eta_0+\eta_1=z+\onm{ad}_{\tilde\xi}\xi_e.
\end{gather*}
As in the proof of Theorem~\rf{normaltransv}, these
conditions are equivalent to the conditions that there exist
$\tilde\xi\in\mfk g$ and $z\in\mfk z_e$ such that
\begin{gather}
\mu_0\in \mfk z_e,\quad C\mu_0=0,\quad dX_{H_{\mu_e}}
(\bar p_e)w=0\lb{eq9}\\
-\onm{ad}_{\xi_e}\xi_0+C^*w+D\mu_0-z-\onm{ad}_{\tilde\xi-\xi_1}
\xi_e=0\lb{eq10}
\end{gather}
subject only to the constraint that $\tilde\xi-\xi_1\in\mfk g_{\mu_e}$,
as long as one takes
\begin{equation}\lb{eq12}
\eta_1=-\onm{ad}_{\xi_e}\xi_0,
\quad\eta_0=z+\onm{ad}_{\tilde\xi-\xi_1}\xi_e.
\end{equation}
Again as in the proof of Theorem~\rf{normaltransv},~\rf{eq9}
and~\rf{eq10} are equivalent to~\rf{eq9} and
\begin{equation}\lb{eq15}
C^*w+D\mu_0-z=0,\quad\onm{ad}_{\tilde\xi-\xi_1-\xi_0}\xi_e=0,
\end{equation}
or equivalently~\rf{eq9}, $C^*w+D\mu_0\in\mfk z_e$ and
$\tilde\xi-\xi_1-\xi_0\in\mfk k_e$. This proves the theorem
since for any $\xi_0$ and $\xi_1$ we can find $\tilde\xi$ such that
$\tilde\xi-\xi_1-\xi_0\in\mfk k_e$.\qed
\end{proof}

One can calculate the linearized dependence at $p_e$ of the generators
of the relative equilibria on the relative equilibria themselves by
taking, at the end of the proof of Theorem~\ref{tanspace},
$\tilde\xi-\xi_1-\xi_0=k\in\mfk k$, and substituting this and $z$
from~\rf{eq15} into~\rf{eq12}, which gives
\begin{equation}\lb{eq16}
\eta_0=C^*w+D\mu_0-\onm{ad}_{\xi_e}\xi_0,\quad
\eta_1=-\onm{ad}_{\xi_e}\xi_1.
\end{equation}
If we define
\begin{equation*}
\hat\mcl E_{(K_e)}=\psi_{(K_e)}^{-1}\bigl(Z(TP)\bigl),
\end{equation*}
that is $\hat\mcl E_{(K_e)}$ is the set of relative equilibria 
paired with their generators, then
\begin{equation*}
T_{p_e}\hat\mcl E_{(K_e)}=\set{(v,\eta)\in T_{p_e}P\times\mfk g}{
dX_{H_{\xi_e}}(p_e)v-\eta.p_e=0},
\end{equation*}
and so $T_{p_e}\hat\mcl E_{(K_e)}$ is the graph of~\rf{eq16} over the tangent
space $T_{p_e}\mcl E_{(K_e)}$ calculated in Theorem~\ref{tanspace}.

The following result uses Theorem~\ref{tanspace} to give a
generalization and partial converse of a result
of~\ct{PatrickGW-1995.1}.

\begin{theorem}\lb{mainth4}
If the relative equilibrium $p_e$ is transversal then
 $\mcl E_{(K_e)}$ is a symplectic
submanifold of $P$ in a neighborhood of $p_e$ if and only if $p_e$ is
nondegenerate and $G_{\mu_e}$ is a maximal torus of $G$.
\end{theorem}

\begin{proof}
The manifold $\mcl E_{(K_e)}$ is symplectic near $p_e$ if and only if
the tangent space $T_{p_e}\mcl E_{(K_e)}$ is a symplectic subspace of
$T_{p_e}P$. From Theorem~\ref{tanspace} and the description of the symplectic
form on $T_{p_e}P$ in terms of the decomposition~\rf{decomp2} it is clear
that the restriction of the symplectic form to $T_{p_e}\mcl E_{(K_e)}$
is degenerate on the subspace
$\onm{ann}_{\mfk g_{\mu_e}}\bigl(\mfk z_e^*\bigr)
\subset \mfk g_{\mu_e} \cong T_0$.
Thus $T_{p_e}\mcl E_{(K_e)}$ can only be a symplectic subspace of $T_{p_e}P$
if $\onm{ann}_{\mfk g_{\mu_e}}\bigl(\mfk z_e^*\bigr)=\{0\}$, that is
$\mfk g_{\mu_e}=\mfk z_e$. This implies that $\mfk g_{\mu_e}=\mfk k_e$ and is
abelian, and so 
$G_{\mu_e}$ is a maximal torus of $G$ by Lemma~\ref{g*gisgs}.

When $\mfk z_e = \mfk k_e = \mfk g_{\mu_e}$ the condition
$C^*w+D\mu_0 \in \mfk z_e$ in the description of $T_{p_e}\mcl E_{(K_e)}$
in Theorem~\ref{tanspace} is automatically satisfied. With
respect to the decomposition~\rf{decomp2} we therefore
have
\begin{equation}\lb{ndtanspace}
T_{p_e}\mcl E_{(K_e)}=\mfk g_{\mu_e}\oplus\ker dX_{H_{\mu_e}}(\bar p_e)
\oplus\ker C\oplus\mfk g_{\mu_e}^\perp.
\end{equation}
The subspace $\mfk g_{\mu_e}^\perp$ is always symplectic. Since $p_e$ is
transversal any $0$ eigenvalue of $dX_{H_{\mu_e}}(\bar p_e)$ is
semisimple and so $\ker dX_{H_{\mu_e}}(\bar p_e)$ is symplectic. It follows
that $T_{p_e}\mcl E_{(K_e)}$ is a symplectic subspace of $T_{p_e}P$ if and
only if $\ker C = \mfk g_{\mu_e}^{*}$ and so $C=0$. However, for a transversal
relative equilibrium $C$ must map onto $\ker dX_{H_{\mu_e}}(\bar p_e)$ and
so $T_{p_e}\mcl E_{(K_e)}$ can only be a symplectic subspace of $T_{p_e}P$
if $p_e$ is nondegenerate.

Conversely, suppose $p_e$ is nondegenerate and $G_{\mu_e}$ is a maximal torus.
Then $K_e\subseteq G_{\mu_e}$, $K_e$ contains a maximal torus and so
$K=G_{\mu_e}$ and $\mfk z_e=\mfk k_e=\mfk g_{\mu_e}$. It follows that
$C=0$, equation~\rf{ndtanspace} holds, and $T_{p_e}\mcl E_{(K_e)}$ is a
symplectic subspace of $T_{p_e}P$. \qed\end{proof}

\subsection*{Acknowledgments}
The authors are grateful to the UK EPSRC for a Visiting Fellowship grant.
G.~W.~Patrick was partially supported by NSERC grant~OGP010571.

\end{document}